\documentclass[a4paper, 11pt]{article}

\usepackage{amsmath}
\usepackage{amsfonts}
\usepackage{epsfig}
\usepackage{graphics} 
\usepackage{authblk}

\newcommand{\ns}{\mathcal H_K(\Omega)}
\newcommand{\R}{\mathbb R}
\newcommand{\N}{\mathbb N}

\usepackage{algorithm,algpseudocode,amsmath}

\newcounter{phase}[algorithm]
\newlength{\phaserulewidth}

\makeatother

\title{Greedy kernel methods for accelerating implicit integrators for parametric ODEs}
\author[1]{T. Br{\"u}nnette}
\author[1]{G. Santin \thanks{santinge@mathematik.uni-stuttgart.de, 
orcid.org/0000-0001-6959-1070}}
\author[1]{B. Haasdonk \thanks{haasdonk@mathematik.uni-stuttgart.de}}

\affil[1]{Institute for Applied Analysis and Numerical Simulation, University of Stuttgart, Germany}

\begin{document}

\maketitle

\begin{abstract}
We present a novel acceleration method for the solution of parametric ODEs by single-step implicit solvers by means of greedy kernel-based surrogate 
models.
In an offline phase, a set of trajectories is precomputed with a high-accuracy ODE solver for a selected set of parameter samples, and used to train a kernel model which 
predicts the next 
point in the 
trajectory 
as a 
function of the last one. This model is cheap to evaluate, and it is used in an online phase for new parameter samples to provide a good initialization point for the 
nonlinear solver of the 
implicit 
integrator. The accuracy of the surrogate reflects into a reduction of the number of iterations until convergence of the solver, thus providing an overall speedup of 
the full simulation.
Interestingly, in addition to providing an acceleration, the accuracy of the solution is maintained, since the ODE solver is still used to guarantee the 
required precision.
Although the method can be applied to a large variety of solvers and different ODEs, we will present in details its use with the Implicit Euler method 
for 
the solution of the Burgers equation, which results to be a meaningful test case to demonstrate the method's features. 
\end{abstract}

\section{Problem setting}\label{santin_mini_MS10_sec:intro}
We consider a $d$-dimensional, autonomous, first order parametric initial value problem: For a given vector of parameters $\mu\in \mathcal P\subset \R^p$ from an 
admissible set $\mathcal P$, solve

\begin{equation*}
\mathbf{IVP(\mu):}\left\{
\begin{array}{ccl}
\dot u(t, \mu) &=& f(u(t,\mu), \mu), \; t\in [0, T]\\
u(0) &=& u_0(\mu)\in\R^d.
\end{array}
\right.
\end{equation*}

We assume that $\mathbf{IVP(\mu)}$ has a unique solution $u(t, \mu):= u(t, \mu, u_0(\mu))$, $t\in [0, T]$, for any value $\mu\in\mathcal P$ and for any initial value  
$u_0(\mu)\in\R^d$. 
Conditions on $f$ such that this requirement is fulfilled are well known, and we refer e.g. to \cite{santin_mini_MS10_Hairer2009} for details.

Existence and uniqueness of solutions allow to define a parametric time evolution or flow mapping 
\begin{equation}\label{santin_mini_MS10_eq:flow}
\Phi(t, u_0(\mu)) := u(t, \mu), 
\end{equation}
which maps the initial value and the time to the corresponding solution vector in $\R^d$, and for which it holds $\Phi(s, u(t, \mu)) = u(t + s, \mu)$.
Although the dependency on the parameters in $\mathbf{IVP(\mu)}$ can be quite general, we require that $\Phi(t, u_0(\mu)) \neq \Phi(s, u_0(\nu))$ for all $t, s\in [0, 
T]$ and for all $\mu, \nu\in\mathcal P$ s.t. $(t,\mu)\neq (s,\nu)$, i.e., different parameters lead to non intersecting trajectories.

We further assume to have an implicit time integrator which is able to numerically solve $\mathbf{IVP(\mu)}$ with any given accuracy, provided a small enough time step 
is 
used.  Although our acceleration algorithm applies to general single-step integration methods, in this paper, for the sake of presentation, we will 
concentrate on the Implicit Euler method (IE), and we refer again to \cite{santin_mini_MS10_Hairer2009} for details on its accuracy.

Such integration method considers a timestep $\Delta t>0$ and a uniform time discretization of $[0, T]$ in $N_t:=N_{\Delta t}\in\N$ intervals $0=t_0<t_1<\dots <t_{N_t} 
\leq T$, with $t_{i+1} - t_{i} = \Delta t$, $0\leq 
i \leq N_t-1$, and computes 
a discrete-time approximation of $u$ as $u_i(\mu)\approx u(t_i, \mu),\; 0\leq i\leq N_t$.

A numerical time evolution map $\phi:\R\times\R^d\to\R^d$ analogous to \eqref{santin_mini_MS10_eq:flow} can be defined from the approximate solution as 
\begin{equation}
\phi(\Delta t, u_i(\mu)) := u_{i + 1}(\mu), \; 0\leq i\leq N_t-1,
\end{equation}
i.e., the solution vector at the current time point is mapped to the solution vector at the next time point. Observe that, under the hypotheses of arbitrary accuracy of 
the integration method and of non intersection of the trajectories, we assume that also the discrete trajectories are 
non intersecting. This means that $\phi$ is a globally defined function independent of the parameter $\mu\in\mathcal P$. 


At each discrete time point, the integrator needs to solve a generally nonlinear, $d$-dimensional system of equations to determine the approximation $u_i(\mu)$. We 
assume 
that this equation is solved with an iterative method, e.g., the Newton method, using an initialization $\bar u_i(\mu)\in\R^d$ at time $t_i$. Common choices 
of this value for the IE method are, e.g., the previous approximation $u_{i-1}(\mu)$ or the approximation obtained by one step of the Explicit Euler method.

The goal of this paper is to present a way to accelerate the computation of the numerical solution $\{u_i(\mu)\}_{i}$ for an arbitrary parameter vector $\mu\in\mathcal 
P$.

The acceleration is realized by constructing a surrogate $s_{\phi}:\R\times\R^d\to\R^d$ of the numerical time evolution map $\phi$ such that $s_{\phi}(\Delta t, u) 
\approx \phi(\Delta t, u)$ for all $(\Delta t, u) \in [0, T]\times \R^d$, while the evaluation of $s_{\phi}$ is much faster than the evaluation of $\phi$.
This surrogate is computed in an offline phase in a data-dependent fashion, i.e., it is trained using a set of precomputed numerical trajectories $\{u_i(\mu_j)\}_{ij}$ 
for multiple parameter values $\mathcal P_{tr}:=\{\mu_1,\dots,\mu_{N_{\mu}}\}\subset\mathcal P$, $N_{\mu}\in\N$ and possibly multiple timesteps $\Delta t$. 

In the online phase, for a new parameter $\mu\in\mathcal P$ the numerical solution is computed by the same time integrator and with timestep $\Delta t$, and, at each 
time $t_i$, the nonlinear 
solver is initialized by $s_{\phi}(\Delta t, u_{i-1})$, i.e., $\bar u_i(\mu)\in\R^d$ is replaced by the surrogate prediction based on the previous timestep.

If the surrogate is accurate, $s_{\phi}(\Delta t, u_{i-1}(\mu))$ is a good approximation of $\phi(\Delta t, u_{i-1}(\mu)) = u_i(\mu)$, so the nonlinear 
solver will converge in possibly significantly less iterations, ideally in $0$ iterations if a residual criterion is used before starting the fix-point loop. This 
reduction of the iterations, combined with the fast evaluation of $s_{\phi}$, will produce a speedup 
of the overall computational time. Moreover, since the same time integrator and nonlinear solver are used in the accelerated algorithm, we should expect no degradation 
of the accuracy, provided the surrogate prediction is accurate enough so that the initialization point is within the area of convergence of the nonlinear 
solver. This is in contrast to general surrogate modeling or model reduction, where the approximation typically results in an accuracy loss.

The surrogate is constructed using the Vectorial Kernel Orthogonal Greedy Algorithm (VKOGA) \cite{santin_mini_MS10_Wirtz2013}, which will be discussed in Section 
\ref{santin_mini_MS10_sec:VKOGA}. In 
particular, it is a kernel-based interpolation algorithm that constructs a nonlinear surrogate $s_{\phi}$.  The full specification of the training data and the complete 
acceleration algorithm will be described in Section \ref{santin_mini_MS10_sec:VKOGA-IE}, but we anticipate that arbitrary unstructured trajectory data 
$\{u_i(\mu_j)\}_{ij}$ in possibly high dimension $d$ can be used. We will conclude this paper with different numerical experiments in Section 
\ref{santin_mini_MS10_sec:numerics} to demonstrate the capabilities of our method.

Moreover, similar acceleration methods have been presented in the papers 
\cite{santin_mini_MS10_CBHB16,santin_mini_MS10_Carlberg201579}, where instead a linear surrogate is employed.

\section{Kernel based surrogates and the VKOGA}\label{santin_mini_MS10_sec:VKOGA}
We briefly outline here the fundamentals of interpolation with kernels and of the VKOGA algorithm, and we refer to \cite{santin_mini_MS10_Wendland2005} and to 
\cite{santin_mini_MS10_Wirtz2013,santin_mini_MS10_HS2017a} 
for the respective details.

We assume to have a function $f:\Omega\subset\R^p\to \R^q$ and a training dataset composed of pairwise distinct data points $X:=\{x_i\}_{i=1}^N\subset\Omega$ 
and data values $Y:=\{f(x_i)\}_{i=1}^N\subset \R^q$. We will specify in the following section the definition of the dataset for the current algorithm.

The general form of the surrogate is
\begin{equation}\label{santin_mini_MS10_eq:ansatz}
 s_{f}(x) :=\sum_{j=1}^N \alpha_j K(x, x_j),\; x\in \Omega,
\end{equation}
where $\alpha_j\in\R^q$ are coefficient vectors and $K:\Omega\times\Omega\to\R$ is a symmetric and strictly positive definite kernel. This means that the matrix $A_{X, 
K}\in\R^{N\times N}$, $\left(A_{X,K}\right)_{ij}:=K(x_i, x_j)$ is positive definite for all $N\in\N$ and for all sets $X\subset\Omega$ of $N$ pairwise distinct points. A 
particular $K$, 
i.e., the Gaussian kernel $K(x, y):=\exp(-\varepsilon ^2 \|x-y\|_2^2)$, with a positive shape parameter
$\varepsilon>0$, will be used in Section \ref{santin_mini_MS10_sec:numerics}.

The coefficient vectors in \eqref{santin_mini_MS10_eq:ansatz} can be 
uniquely 
determined by imposing interpolation conditions 
\begin{equation}
 s_{f}(x_i) :=f(x_i), \;\; 1\leq i\leq N,
\end{equation}
which result, defining $\alpha^T := \left[\alpha_1, \dots, \alpha_N\right]$, $b^T := \left[f(x_1), \dots,f(x_N)\right]$, $\alpha, b \in\R^{q\times N}$, 
in the solution of the linear system $A_{X, K}\ \alpha = b$. This, indeed, has a unique solution as $A_{X, K}$ is positive definite by assumption.

This interpolation method is well studied, and we just recall that convergence rates are proven for functions $f$ in the space $\ns$, which is a Reproducing 
Kernel Hilbert Space (RKHS) associated to the particular kernel $K$, and which is norm equivalent to a Sobolev space $W_2^{\tau}(\Omega)$, $\tau>d/2$, for certain 
kernels (see \cite{santin_mini_MS10_Wendland2005}).

The goal of the VKOGA is to approximate the surrogate \eqref{santin_mini_MS10_eq:ansatz} by a sparse expansion of the same form, i.e., one where most of the $\alpha_j$ 
are the zero 
vector. A good selection of the sparsity pattern results into an approximate surrogate which is as good as the full one, while being much faster to evaluate, since the 
sum involves only $n\ll N$ elements. The selection of the non-zero coefficients and their computation is realized by a greedy procedure in $\ns$, which 
iteratively selects nested data point sets $\emptyset\subset X_{n-1}\subset X_{n}\subset\Omega$ by maximizing a selection criterion at each step, and solves the 
corresponding interpolation problem. Possible choices in the VKOGA are the $f$-, $P$-, and $f/P$-greedy selection rules 
\cite{santin_mini_MS10_SchWen2000,santin_mini_MS10_DeMarchi2005,santin_mini_MS10_Muller2009}. The 
algorithm has theoretical grounds, e.g. provable convergence rates \cite{santin_mini_MS10_Muller2009,santin_mini_MS10_Wirtz2013}, which are also quasi-optimal in Sobolev 
spaces for $P$-greedy 
\cite{santin_mini_MS10_SH16b}, and has been successfully applied in several application contexts, e.g. \cite{santin_mini_MS10_KSHH2017}. Moreover, the numerical 
computation of the 
surrogate can be 
efficiently implemented using a partial Cholesky decomposition of the kernel matrix $A_{X, K}$, where only the columns appearing in the sparse surrogate need to be 
computed and stored.

\section{The complete algorithm: VKOGA-IE}\label{santin_mini_MS10_sec:VKOGA-IE}
We can now describe the complete algorithm, which we name VKOGA-IE.
The target function is $f:=\phi$, which is defined on $\Omega:=[0,T]\times \R^d$ to $\R^d$, i.e., $p := d+1$, $q := d$. 
What remains to specify is the exact definition of the training set $(X, Y)$ used by VKOGA to construct the surrogate $s_{\phi}$, as described in the previous section. 
As mentioned in Section \ref{santin_mini_MS10_sec:intro}, we solve $\mathbf{IVP}(\mu)$ for $N_{\mu}\in\N$ different parameters from a parameter-training set $\mathcal 
P_{tr}\subset 
\mathcal P$, each $\mu_j$ with a timestep $\Delta t_j$. If the same parameter is used more than once with different timesteps, we just count it multiple times in 
$\mathcal P_{tr}$. This generates trajectory data which we assign at temporary sets $X_j:= \{(\Delta t_j, u_i(\mu_j))\}_{i=0}^{N_t-1}$, $Y_j:= 
\{u_{i+1}(\mu_j)\}_{i=0}^{N_t-1}$, representing input-output pairs of $\phi$. The dataset is defined as $X:=\cup_{j=1}^{N_{\mu}} X_j$, $Y:=\cup_{j=1}^{N_{\mu}} Y_j$. 
The complete offline phase is summarized in Algorithm \ref{santin_mini_MS10_alg:offline}.

Instead of working with a fixed kernel shape parameter $\varepsilon>0$, typically step $13$ implies a parameter selection procedure, e.g. via cross validation.
Moreover, we assume for simplicity that $T/\Delta t\in\N$.
\begin{algorithm}[ht!]
  \caption{VKOGA-IE (Offline phase)}
  \begin{algorithmic}[1]
      \State{Input: $\{\Delta t_j, \mu_j\}_{j=1}^{N_{\mu}}$}
      \For{$i = 1, \dots, N_{\mu}$}
        \State{$u_0(\mu_j):= u(0, \mu_j)$}
        \For{$i = 1, \dots, N_{t}$}
	\State{Initialize $\bar u_i(\mu)$}
	\State{Compute $u_i(\mu)$ with IE}
        \EndFor
        \State{$X_j:= \{(\Delta t_j, u_i(\mu_j))\}_{i=0}^{N_t-1}$}
        \State{$Y_j:= \{u_{i+1}(\mu_j)\}_{i=0}^{N_t-1}$}
      \EndFor
      \State{$X:=\cup_{j=1}^{N_{\mu}} X_j$}
      \State{$Y:=\cup_{j=1}^{N_{\mu}} Y_j$}
      \State{Train $s_{\phi}$ on dataset $(X, Y)$ with VKOGA}
      \State{Output: $s_{\phi}$}
    \end{algorithmic}\label{santin_mini_MS10_alg:offline}
\end{algorithm}
In the online phase, instead, we only need to run the IE method and solve at each iteration the nonlinear equation using the initialization provided by the surrogate, as 
described in Algorithm \ref{santin_mini_MS10_alg:online}.
\begin{algorithm}[ht!]
  \caption{VKOGA-IE (Online phase)}
  \begin{algorithmic}[1]
      \State{Input: $\Delta t$, $\mu$, $s_{\phi}$}
      \State{$\mu\in\mathcal P$,  $\Delta t>0$}
        \State{$u_0(\mu):= u(0, \mu)$}
        \For{$i = 1, \dots, N_{t}$}
	\State{Initialize $\bar u_i(\mu)=s_{\phi}(\Delta t, u_{i-1}(\mu))$}
	\State{Compute $u_i(\mu)$ with IE}
        \EndFor
     \State{Output: $\{u_i(\mu)\}_{i=0}^{N_t}$}
    \end{algorithmic}\label{santin_mini_MS10_alg:online}
\end{algorithm}

\section{Experiments}\label{santin_mini_MS10_sec:numerics}
To demonstrate the features of VKOGA-IE, we consider the Burgers equation 
\begin{equation*}
\left\{\begin{array}{lcll}
\partial_t \theta(t, x) + \frac{1}{2}\partial_x \theta(t, x) ^2 &=& 0,& (t, x) \in [0, T] \times [-r, r]\\
\theta(0,x) &=& \theta_0(x)&x\in [-r, r]\\
\theta(t,-r) &= &u_l,&t\in [0, T]\\
\theta(t,r) &=& u_r,& t\in [0, T],
\end{array}
\right.
\end{equation*}
which is transformed into an ODE by a semi-discrete finite volume discretization in space based on the Lax-Friedrichs flux. We consider $d:=200$ cells in $(-r, r)$ with 
$r:=5$, and $T$ as specified later. This produces a $d$ dimensional $\mathbf{ IVP(\mu)}$ depending on a two-dimensional parameter vector $\mu:=(u_l, u_r)$. We 
concentrate here on shock 
wave solutions, i.e., 
$u_l>u_r$. The resulting ODE is then simulated from $t=0$ to $t=T$, with varying time-step $\Delta t$. The nonlinear system is solved at each timestep using the Newton 
method, which is 
terminated with a maximal number of $100$ iterations or when a tolerance of $10^{-14}$ on the residual is reached. To have more training points, all the training sets 
in the following are generated with training time $T = T_{tr}:=4$. 

The VKOGA is run with the Gaussian kernel and with a termination tolerance of $10^{-12}$. The kernel depends on a parameter $\varepsilon>0$, which is chosen via 
$5$-fold cross validation from a set of $50$ logarithmically equally spaced values in $[10^{-4}, 10^{2}]$. 

The first experiment uses a fixed $\Delta t=0.01$ and a single training parameter $\mathcal P_{tr} = \{(3.4, 0.2)\}$, i.e., $N = 400 = T/\Delta t$. Observe 
that a fixed $\Delta t$ means that the model is in practice $d$ to $d$ dimensional. The VKOGA selects $n = 67$ points, and the model is tested to solve $\mathbf{ 
IVP(\mu)}$ with parameters $\mathcal P_{te}:=\{(3.4 + i\ 0.2, 0.2 + j\ 0.2), i, j\in\{-1,0, 1\}\}$ and $T = T_{te}:=2$.
The results are summarized in Table \ref{santin_mini_MS10_tab:reproduction}. The average number of iterations for the standard initialization with the previous value 
('Old value' column) 
and with the VKOGA model ('VKOGA' column) are reported, as well as the test parameters where the minimal and maximal gain of our technique is realized. The table 
contains also the computational times in seconds, which are the averages over $10$ repetitions of the same simulation based on a Matlab implementation.
It is evident that a good speedup is reached when the model is tested on the training parameter, while the quality degrades for different ones.

\begin{table}[ht!]
\begin{center}
\begin{tabular}{||c||c|c||c|c||c|c||c||}
\hline
& \multicolumn{2}{|c||}{Old value}& \multicolumn{2}{|c||}{VKOGA}&\multicolumn{2}{|c||}{Gain}&\\
\hline
&iter&time&iter&time&iter&time&$\mu$\\
\hline
Mean&$25.40$ & $1.38$ &  $25.83$&$1.50$ & $-1.70\%$&$-8.55\%$&\\
\hline
Min&$24.27$&$1.29$&$27.84$&$1.63$&$-14.69\%$&$-25.97\%$&$(3.2, 0.4)$\\
\hline
Max&$25.30$ &$1.34$&$16.41$&$0.98$&$35.11\%$&$26.59\%$&$(3.4, 0.2)$\\
\hline
\end{tabular}
\end{center}
\caption{Results of the first experiment with $\mathcal P_{tr} = \{(3.4, 0.2)\}$ and fixed timestep $\Delta t = 0.01$.}\label{santin_mini_MS10_tab:reproduction}
\end{table}

\vspace*{-.5cm}The second experiment uses a model trained again with fixed $\Delta t=0.01$ and the same test parameters $\mathcal P_{te}$, but instead with training 
parameters 
$\mathcal P_{tr}:=\{(3.2, 0), (3.2, 0.4), (3.6, 0), (3.6, 0.4)\}$, i.e., the corners of the square containing $\mathcal P_{te}$. The resulting training set has $N = 400 
\times 4 = 1600$ points, and the VKOGA selects $n= 219$ points. The results are summarized in Table \ref{santin_mini_MS10_tab:generalization}. In this case, as expected, 
we obtain a 
significant reduction of the number of iterations for all the test parameters. The minimal reduction is realized for the parameter $\mu=(3.2, 0.2)$, which is the 
farthest from the training set. This is a further indication that the quality of the model degrades with the distance from the training set, which is a reasonable 
behavior 
but also a promising feature, since a model trained on a larger parameter training set should improve the acceleration. This reduction is reflected also in a speedup in 
terms of computational cost, except in one case reported in the table. This suggests that the additional cost required by the evaluation of the kernel model is relevant 
in the case of a small reduction of the number of iterations. Nevertheless, the computational time is highly dependent on the implementation, while the number of 
iterations is not.

\begin{table}[ht!]
\begin{center}
\begin{tabular}{||c||c|c||c|c||c|c||c||}
\hline
& \multicolumn{2}{|c||}{Old value}& \multicolumn{2}{|c||}{VKOGA}&\multicolumn{2}{|c||}{Gain}&\\
\hline
&iter&time&iter&time&iter&time&$\mu$\\
\hline
Mean&$25.40$ & $1.36$ &  $20.31$&$1.17$ & $19.98\%$&$13.35\%$&\\
\hline
Min&$25.29$&$1.28$&$24.25$&$1.33$&$4.13\%$&$-4.18\%$&$(3.4, 0.2)$\\
\hline
Max&$25.12$ &$1.36$&$16.68$&$0.97$&$33.60\%$&$28.21\%$&$(3.2, 0)$\\
\hline
\end{tabular}
\end{center}
\caption{Results of the second experiment, i.e., model trained with $\mathcal P_{tr}:=\{(3.2, 0), (3.2, 0.4), (3.6, 0), (3.6, 
0.4)\}$ and fixed timestep $\Delta t = 0.01$.}\label{santin_mini_MS10_tab:generalization}
\end{table}

\vspace*{-.5cm}Finally, we test the behavior of the method with respect to a change in the timestep $\Delta t$. To this end, we use $\mathcal P_{tr} = \mathcal P_{te} = 
\{(3.4, 
0.2)\}$, but we train the model with the solutions computed for $\Delta t \in\{0.01, 0.005, 0.001\}$ and test for $10$ logarithmically equally spaced timesteps in 
$[0.001, 0.05]$. The results are reported in Table \ref{santin_mini_MS10_tab:variable}. Also in this case we achieve a reduction of the number of Newton iterations in all 
cases, even if 
this reduction is not sufficient in the case of the smallest timestep to achieve a computational speedup, since the number of iterations is already quite small. 
Nevertheless, the reduction of the number of iterations suggests that the kernel model captures well the dependence on the timestep, so one could expect to use this 
technique in more general settings without the need of including in the training sets many solutions obtained with different timesteps.
\begin{table}[ht!]
\begin{center}
\begin{tabular}{||c||c|c||c|c||c|c||c||}
\hline
& \multicolumn{2}{|c||}{Old value}& \multicolumn{2}{|c||}{VKOGA}&\multicolumn{2}{|c||}{Gain}&\\
\hline
&iter&time&iter&time&iter&time&$\Delta t$\\
\hline
Mean&$32.45$ & $2.28$ &  $29.20$&$2.29$ & $11.29\%$&$4.22\%$&\\
\hline
Min&$9.20$&$5.74$&$8.98$&$6.46$&$2.45\%$&$-12.44\%$&$10^{-3}$\\
\hline
Max&$23.98$ &$1.07$&$19.21$&$0.94$&$19.88\%$&$12.51\%$&$8.79 10 ^{-3}$\\
\hline
\end{tabular}
\end{center}
\caption{Results of the third experiment, i.e., $\mathcal P_{tr} = \{(3.4, 0.2)\}$ and multiple timesteps.}\label{santin_mini_MS10_tab:variable}
\end{table}

\vspace*{-1cm}\section{Conclusion and further work}
In this work we described a general nonlinear forecasting method used for the acceleration of implicit ODE integrators. The method is suited for parametric problems and 
multi-query 
scenarios, and it realizes a significant acceleration possibly without accuracy degradation. 

The algorithm can be extended to non-autonomous ODEs, adaptive-timestep or multi-stage Runge-Kutta time integrators. In each case, more simulation data should be 
included in the training set, such as the current time or the partial solutions of the intermediate stages.

Another interesting aspect that could be investigated is the analysis of the accuracy of the method. Indeed, if it is possible to prove that the surrogate has a small 
enough uniform 
error, it would be guaranteed that the initialization point is inside the convergence area of the nonlinear solver.

\section*{Acknowledgments}  The authors would like to thank the German Research Foundation (DFG) for financial support of the project within the Cluster of Excellence in
Simulation Technology (EXC 310/2) at the University of Stuttgart.


\begin{thebibliography}{99}
\parskip1.0ex

\bibitem{santin_mini_MS10_CBHB16}
{\sc K.~Carlberg, L.~Brencher, B.~Haasdonk, and A.~Barth}, {\em Data-driven
  time parallelism via forecasting}, submitted to SIAM J. of Sci. Comp., 2016.

\bibitem{santin_mini_MS10_Carlberg201579}
{\sc K.~Carlberg, J.~Ray, and B.~van Bloemen~Waanders}, {\em Decreasing the
  temporal complexity for nonlinear, implicit reduced-order models by
  forecasting}, Computer Methods in Applied Mechanics and Engineering {\bf 289}
  (2015), 79 -- 103.

\bibitem{santin_mini_MS10_DeMarchi2005}
{\sc S.~De~Marchi, R.~Schaback, and H.~Wendland}, {\em Near-optimal
  data-independent point locations for radial basis function interpolation},
  Adv. Comput. Math. {\bf 23}:3 (2005), 317--330.

\bibitem{santin_mini_MS10_HS2017a}
{\sc B.~Haasdonk and G.~Santin}, {\em Greedy kernel approximation for sparse
  surrogate modelling}, Proceedings of the KoMSO Challenge Workshop on
  Reduced-Order Modeling for Simulation and Optimization, 2017.

\bibitem{santin_mini_MS10_Hairer2009}
{\sc E.~Hairer, S.~P. N\o~rsett, and G.~Wanner}, {\em Solving {O}rdinary
  {D}ifferential {E}quations. {I}: Nonstiff problems.}, second ed., Springer Series in
  Computational Mathematics, vol.~8, Springer-Verlag, Berlin, 1993, 

\bibitem{santin_mini_MS10_KSHH2017}
{\sc T.~K{\"o}ppl, G.~Santin, B.~Haasdonk, and R.~Helmig}, {\em Numerical
  modelling of a peripheral arterial stenosis using dimensionally reduced
  models and machine learning techniques}, Tech. report, University of
  Stuttgart, 2017.

\bibitem{santin_mini_MS10_Muller2009}
{\sc S.~M{\"u}ller and R.~Schaback}, {\em A {N}ewton basis for kernel spaces},
  J. Approx. Theory {\bf 161}:2 (2009), 645--655.

\bibitem{santin_mini_MS10_SH16b}
{\sc G.~Santin and B.~Haasdonk}, {\em Convergence rate of the data-independent
  {P}-greedy algorithm in kernel-based approximation}, Dolomites Res.\ Notes Approx. 
  {\bf 10} (2017), 68--78.

\bibitem{santin_mini_MS10_SchWen2000}
{\sc R.~Schaback and H.~Wendland}, {\em Adaptive greedy techniques for
  approximate solution of large {RBF} systems}, Numer. Algorithms {\bf 24}:3
  (2000), 239--254.

\bibitem{santin_mini_MS10_Wendland2005}
{\sc H.~Wendland}, {\em Scattered {D}ata {A}pproximation}, Cambridge Monographs
  on Applied and Computational Mathematics, vol.~17, Cambridge University
  Press, Cambridge, 2005.

\bibitem{santin_mini_MS10_Wirtz2013}
{\sc D.~Wirtz and B.~Haasdonk}, {\em A vectorial kernel orthogonal greedy
  algorithm}, Dolomites Res. Notes Approx. {\bf 6} (2013), 83--100.

\end{thebibliography}
\end{document}